\documentclass[12pt]{amsart}
\usepackage{amssymb,latexsym,comment,url}

\newcommand{\PP}{{\mathbb P}}
\newcommand{\Q}{{\mathbb Q}}

\newenvironment{Proof}{\par\noindent{\sc Proof:}}%
                      {\hspace*{\fill}\nobreak$\Box$\par\medskip}
                       {\hspace*{\fill}\nobreak$\Box$\par\medskip}

\newtheorem{Proposition}{Proposition}[section]
\newtheorem{Theorem}[Proposition]{Theorem}

\newtheorem{Conjecture}[Proposition]{Conjecture}

\theoremstyle{definition}

\newtheorem{Remark}[Proposition]{Remark}

\renewcommand{\baselinestretch}{1.1}

\begin{document}

\title[Polynomials with no preperiodic points]
{Families of polynomials of every degree with no rational preperiodic points\\
Familles de polyn\^{o}mes de degr\'{e} arbitraire sans points pr\'{e}p\'{e}riodiques rationnels}%

\author[M. Sadek]%
{Mohammad~Sadek}
\address{Faculty of Engineering and Natural Sciences,
 Sabanc{\i} University,
  Tuzla, \.{I}stanbul, 34956 Turkey}
\email{mohammad.sadek@sabanciuniv.edu}
\let\thefootnote\relax
\begin{abstract}
Let $K$ be a number field. Given a polynomial $f(x)\in K[x]$ of degree $d\ge 2$, it is conjectured that the number of preperiodic points of $f$ is bounded by a uniform bound that depends only on $d$ and $[K:\Q]$. However, the only examples of parametric families of polynomials with no preperiodic points are known when $d$ is divisible by either $2$ or $3$ and $K=\Q$. In this article, given any integer $d\ge 2$, we display infinitely many parametric families of polynomials of the form $f_t(x)=x^d+c(t)$, $c(t)\in K(t)$, with no rational preperiodic points for any $t\in K$.
\end{abstract}

\maketitle

{\bf R\'{e}sum\'{e}.}
Soit $K$ un corps de nombres. \'{E}tant donn\'{e} un polyn\^{o}me $f(x)\in K [x]$ de degr\'{e} $d\ge 2$, il est conjectur\'{e} que le nombre de points pr\'{e}p\'{e}riodiques de $f$ est born\'{e} par une constante ne d\'{e}pendant que de $d$ et $[K:\Q]$. Cependant, les seuls exemples de familles param\'{e}triques de polyn\^{o}mes sans points pr\'{e}p\'{e}riodiques supposent $2|d$ ou $3|d$ et $K=\Q$. Dans cet article, \'{e}tant donn\'{e} un entier $d\ge 2$,
nous d\'{e}montrons qu'il existe une infinit\'{e} de familles param\'{e}triques de polyn\^{o}mes de la forme $f_t (x) = x^d + c (t)$, $c (t)\in K(t)$, sans points pr\'{e}p\'{e}riodiques rationnels pour tout $t\in K$.

\section{Introduction}
An arithmetic dynamical system over a number field $K$ consists of a rational function $f: \PP^n(K)\to \PP^n(K)$ of degree at least $2$ with coefficients in $K$ where the $m^{th}$ iterate of $f$ is defined recursively by $f^{1}(x)=f(x)$ and $f^{m}(x)=f(f^{m-1}(x))$ when $m\ge 2$. A point $P\in\PP^n( K)$ is said to be a {\em periodic} point for $f$ if there exists a positive integer $m$ such that $f^m(P)=P$. If $N$ is the smallest positive integer such that $f^N(P)=P$, then the periodic point $P$ is said to be of {\em exact period} $N$. A point $P\in\PP^n( K)$ is said to be a {\em preperiodic} point for $f$ if the orbit $\{f^i(P):i=0,1,2,\cdots\}$ of $P$ is finite, i.e., if some iterate $f^i(P)$ is periodic.

The following conjecture was proposed by Morton and Silverman in p.~4 of \cite{Morton}.
\begin{Conjecture}There exists a bound $B(D, n, d)$ such that if $K/\Q$ is a number field of degree $D$, and $f : \PP^n(K) \to \PP^n(K)$ is a morphism of degree $d\ge 2$ defined over $K$, then the number of $K$-rational preperiodic points of $f$ is bounded by $B(D,n,d)$.
\end{Conjecture}
When $f$ is taken to be a quadratic polynomial over $\Q$, the following conjecture was suggested in \cite[Conjecture 2]{Flynn} and \cite[Conjecture 2]{Poonen}.
\begin{Conjecture}If $N\ge 4$, then there is no quadratic polynomial $f(x)\in\Q[x]$ with a
rational point of exact period $N$.
 \end{Conjecture}
 The conjecture has been proved when $N=4$, see \cite[Theorem 4]{Morton2}, and $N=5$, see \cite[Theorem 1]{Flynn}. A conditional proof for the case $N=6$ was given in \cite[Theorem 7]{Stoll2}.

 Although polynomials described by an equation of the form $x^d+c,\,d\ge2,\,c\in K,$ with rational preperiodic points are scarce, examples of parametric families of polynomials with no preperiodic points are very few in the literature. In Theorem 4 of \cite{Ingram}, families of such polynomials were given when $d$ is even or when $d$ is divisible by $3$, and $K=\Q$. The main finding of this article can be described as follows. Let $K$ be a number field. Given an arbitrary integer $d\ge 2$, we prove the existence of infinitely many parametric families of polynomials of degree $d$ with no $K$-rational preperiodic points. This is achieved using some recent results on the non existence of rational points on certain twisted superelliptic curves.

\subsection*{Acknowledgments} The author would like to thank the referee for several corrections, comments and valuable suggestions that improved the manuscript. He would also like to thank W{\l}adys{\l}aw Narkiewicz for pointing out the papers \cite{Narkiewicz} and \cite{Narkiewicz2}. The author is supported by the starting project B.A.CF-19-01964, Sabanc{\i} University.

\section{Parametric families of polynomials with no periodic points}
In what follows $K$ will denote a number field with ring of integers $\mathcal{O}_K$.
The following proposition is \cite[Lemma 1]{Narkiewicz2}.
\begin{Proposition}
\label{prop:denominators}
Let $f(x)=x^d+c$, where $d\ge 2$ is an integer and $c\in K\setminus\{0\}$. If
 $f$ has a $K$-rational periodic point, then there exist $a,b\in\mathcal{O}_K$ such that $c=a/b^d$, and $(a\mathcal{O}_K,b^d\mathcal{O}_K)=I^d$ for some ideal $I$ in $\mathcal{O}_K$.
\end{Proposition}

  Proposition \ref{prop:denominators} shows that polynomials described by an equation of the form $x^d+c,\,d\ge2,\,c\in K,$ with $K$-rational preperiodic points are rare. However, up to the knowledge of the author, the only such family is given as Theorem 4 in \cite{Ingram} where $K=\Q$. The statement of the latter theorem is as follows.

\begin{Theorem}
Let $2 \mid d$ and $m\ge  4$; or $3 \mid d$ and $m\ge 3$. Then for $t\in\Q$, the polynomial \[x^d +\frac{1}{1 + t^m}\]
has no $\Q$-rational preperiodic points.
\end{Theorem}

Now we state the main result of this work.

\begin{Theorem}
\label{Thm1}
Let $K$ be a number field with ring of integers $\mathcal{O}_K$. Let $d\ge 2$ be an integer. Let $P(T)\in \mathcal{O}_K[T]$ be of degree $N$ a multiple of $d$ such that the multiplicity of each of its roots is
at most $d-1$. Assume moreover that the Galois group of $P(T)$ over $K$ has an element fixing no root of
$P(T)$. Then there exists $w\in\mathcal{O}_K\setminus\{0\}$ such that the polynomial
\[f(x)=x^d+\frac{1}{w\cdot P(t)}\] has no $K$-rational preperiodic points for any $t\in K$.
\end{Theorem}
\begin{Proof}
According to \cite[Theorem 3.1]{Legrand}, given that the Galois group of $P(T)$ over $K$ has an element fixing no root of $P(T)$, it follows that there exists $w\in\mathcal{O}_K\setminus\{0\}$ such that the twisted superelliptic curve defined by $y^d=w\cdot P(T)$ has no $K$-rational points. In other words, there exists no $(y,t,s)\in K^3\setminus\{(0,0,0)\}$ such that $y^d=w\cdot Q(t,s)$, where $Q(T,S)=S^N\cdot P(T/S)$. In view of Proposition \ref{prop:denominators}, $f$ has no $K$-rational periodic points of any period, hence no $K$-rational preperiodic points for any $t\in K$.
\end{Proof}

\begin{Remark}
   One knows that the proportion of degree $N$ polynomials $P(T )\in\mathcal{O}_K[T]$ with height bounded by $H$ and such that the Galois group of $P(T)$ over $K$ is isomorphic to the symmetric group $S_N$ tends to $1$ as $H$ tends to $\infty$, see for example \cite[Theorem 2.1]{Cohen}. Consequently, the proportion of fixed degree polynomials $P(T)$ introduced in Theorem \ref{Thm1} with height bounded by $H$ tends to $1$ as $H$ tends to $\infty$.
\end{Remark}

\end{document}